\documentclass[12pt]{article}

\usepackage{latexsym,amsthm,amsfonts,amsmath,amssymb,amscd} % Математические дополнения от AMS

\usepackage[dvips]{graphicx}

\usepackage[english]{babel} % Языки: русский, английский
\usepackage{textcase}                % Формирование больших букв

\theoremstyle{plain}
\newtheorem{theorem}{Theorem}[section]
\newtheorem{lemma}{Lemma}[section]

\usepackage{indentfirst}

\def \limsup{\mathop{\overline {\rm lim}}}

\begin{document}

\begin{center}
\textbf {Non-local elasticity theory as a continuous limit of 3D
networks of pointwise interacting masses}

         E. Khruslov, M. Goncharenko
\end{center}

\emph{ Small oscillations of an elastic system of point masses
(particles) with a nonlocal interaction are considered. We study
the asymptotic behavior of the system, when number of particles
tends to infinity, and the distances between them and the forces
of interaction tends to zero. The first term of the asymptotic is
described by the homogenized system of equations, which is a
nonlocal model of oscillations of elastic medium.}

\vskip1cm

\begin{Large}
{\bf Introduction}
\end{Large}

The progress in development of new materials and the modelling of
nanostructures caused the emergence of nonlinear elasticity
theories (see, for example, \cite{1}, \cite{2}, \cite{3} ).
Classical local theory is based on the concept of contact
interaction and it can not explain some observed experimental
phenomena. Therefore, it is necessary to take into account the
long-range interaction between the particles of the material and
this leads to the nonlocal elasticity theory.

The nonlocal elasticity theory can be traced back to the works of
Kr\"oner, who formulated the continual theory of elastic materials
with long-range interaction forces (\cite{4}, \cite{5}). At
present, the nonlocal mechanics of the elastic continuum is
treated with two different approaches: the gradient elasticity
theory (weak nonlocality) and the integral nonlocal theory (strong
nonlocality).

The first approach is related to the study of the gradients of the
strain tensors. It leads to models with spatial derivatives of
order more than 2 (\cite{6} - \cite{8}). The main difficulties in
using this model are the setting of boundary conditions for the
corresponding boundary value problems (see \cite{9}).

The second approach has been developed almost independently. The
nonlocal interaction here is represented in the form of a
convolution integral of the deformation tensor with a kernel that
depends on the distance between the particles of the elastic
material. This approach leads to models described by
integro-differential equations (\cite{10} - \cite{13}).

The correctness of these continuum models of nonlocal elasticity
theory depends on the effectiveness of long-range molecular forces
in the material. Therefore, a natural approach to their
justification is the so-called microstructural approach, which is
studying discrete elastic systems (lattice models). This approach
has been used mainly in physical works (\cite{14} - \cite{18}).
Apparently, one of the first mathematical works, in which the
system of equations of the local elasticity theory was derived
using the microstructural approach, was \cite{19}. The short-range
interactions between particles were considered. Only the nearest
particles interact in the system. The asymptotic behavior of the
oscillations of such a system was investigated when the distances
between the nearest neighbors and the forces of interaction
between them tend to zero. A homogenized system of differential
equations describing the leading term of the asymptotic was
obtained. This system is a continuum model of the local theory of
elasticity. In this work the method based on the studying of the
asymptotic behavior of the system, when the scale of the
microstructure tends to zero, was applied. This approach is the
basis for the homogenization of partial differential equations
(\cite{20} - \cite{22}).

We apply this approach of homogenization to study the asymptotic
behavior of the oscillations of an elastic system of point masses
(particles) with a nonlocal interaction. It is assumed that the
system depends on the small parameter $\varepsilon $. More
precisely, the distance between the nearest neighbors is of the
order $O (\varepsilon) $, and the long-range forces are of order
$O (\varepsilon^6) $. It is proved that the main term of the
asymptotic is described by a homogenized system of
integro-differential equations. The integral term is a convolution
of the difference of the displacements of the elastic medium at
various points with some kernel. Note that such a system differs
from the continual model of Eringen, where the convolution of the
deformation tensor with the kernel is taken. A similar system of
integro-differential equations was proposed earlier (without
justification) in \cite{23} as a variant of the integral
elasticity theory and was used to calculate steel plates. The
indicated order of interactions in the system corresponds to the
integral (and not gradient) elasticity theory.

\section{Statement of the problem}

We consider a system $M_\varepsilon $ of interacting point masses
(we will call them particles) in a fixed bounded domain $\Omega
\subset \mathbb{R}^3 $ with a smooth boundary $\partial \Omega $.
It is assumed that this system depends on the small parameter
$\varepsilon > 0 $. The total number of particles in the system is
$ O (\varepsilon^{-3 }) $ and the distances between the nearest
particles are of order $ O(\varepsilon) $. We denote by $
x^i_\varepsilon $ ($ i = 1, ..., N_\varepsilon $) the positions of
the particles in the equilibrium state of the system $
M_\varepsilon $, and we denote by $ u^i_\varepsilon =
u^i_\varepsilon (t) $ the displacements of particles relative to
their equilibrium positions $ x^i_\varepsilon $.

The potential energy for small variation of the system
$M_\varepsilon $ from the equilibrium position is determined by
the equality
\begin{equation}
 H_\varepsilon (u_\varepsilon) = H_0 +\frac{1}{2}
\sum_{i,j=1}^{N_\varepsilon}  \langle E_\varepsilon^{ij}
(u_\varepsilon^i - u_\varepsilon^j), (u_\varepsilon^i -
u_\varepsilon^j)\rangle, \quad H_0={\text{const}}, \label{1.1}
\end{equation}
where $u_\varepsilon = \{
u_\varepsilon^1,...,u_\varepsilon^{N_\varepsilon}\}$,  parentheses
$\langle,\rangle$ denote the scalar product in $\mathbb{R}^3$, and
$E_\varepsilon^{ij}$ are symmetric nonnegative matrices of the
pair interaction between the $i$-th and $j$-th particles. If the
particles interact through the central elastic forces (for
example, they connected by elastic springs), then the matrices
$E_\varepsilon^{ij} $ satisfy the equalities
\begin{equation}
E_\varepsilon^{ij} u=K_\varepsilon^{ij} \langle
u,e_\varepsilon^{ij}\rangle e^{ij}_\varepsilon,\quad \forall
u\in\mathbb{R}^3, \label{1.2}
\end{equation}
where $e^{ij}_\varepsilon =(x_\varepsilon^i -x_\varepsilon^j)\vert
x_\varepsilon^i -x_\varepsilon^j\vert^{-1}$ is the unit vector of
direction between the $i$-th and $j$-th particles and the
coefficient $K_\varepsilon^{ij} $ characterizes the intensity of
interaction (stiffness of springs).

The coefficient $K^{ij}_\varepsilon$ depends on the distances
$\vert x_\varepsilon^i -x_\varepsilon^j \vert $ between particles.
Generally speaking, it can be zero if the corresponding pair of
particles does not interact with each other. In this paper we
assume that the coefficient $ K_\varepsilon^{ij}$ is defined by
formula
\begin{equation}
K_\varepsilon^{ij}=\varepsilon^6 \left [ K(\vert x_\varepsilon^i
-x_\varepsilon^j\vert)  +\frac{K^{ij}}{\vert x_\varepsilon^i
-x_\varepsilon^j\vert^5}\varphi\left (\frac{\vert x_\varepsilon^i
-x_\varepsilon^j\vert}{\varepsilon}\right )\right ]
A_\varepsilon^{ij}, \label{1.3}
\end{equation}
where $K(r)$, $\varphi(r) \in C([0,L])$, $K(r) \ge 0$,
$\varphi(r)=1$ as $r\le \alpha$ and $\varphi(r)=0$  as $r\ge
\beta$ ($0<\alpha<\beta <L={\text{diam}}\, \Omega$);
$A_\varepsilon^{ij}=1$  (for interacting pairs of particles) and
$A_\varepsilon^{ij}=0$ (for noninteracting pairs of particles),
$a_0 \le K^{ij}\le A_0$.

The formula above simulates a weak interaction (of the order $O
(\varepsilon^6)$) between not very close particles ($\vert
x_\varepsilon^i - x_\varepsilon^j\vert >\beta \varepsilon$)  and
stronger interaction $ O (\varepsilon) $ between close ones
($\vert x_\varepsilon^i - x_\varepsilon^j\vert <
\alpha\varepsilon$)  (see Figure 1). This type of interaction is
characteristic for some intermolecular forces (for example, van
der Waals forces).

\begin{figure}[h]
\center{\includegraphics[width=0.5\linewidth]{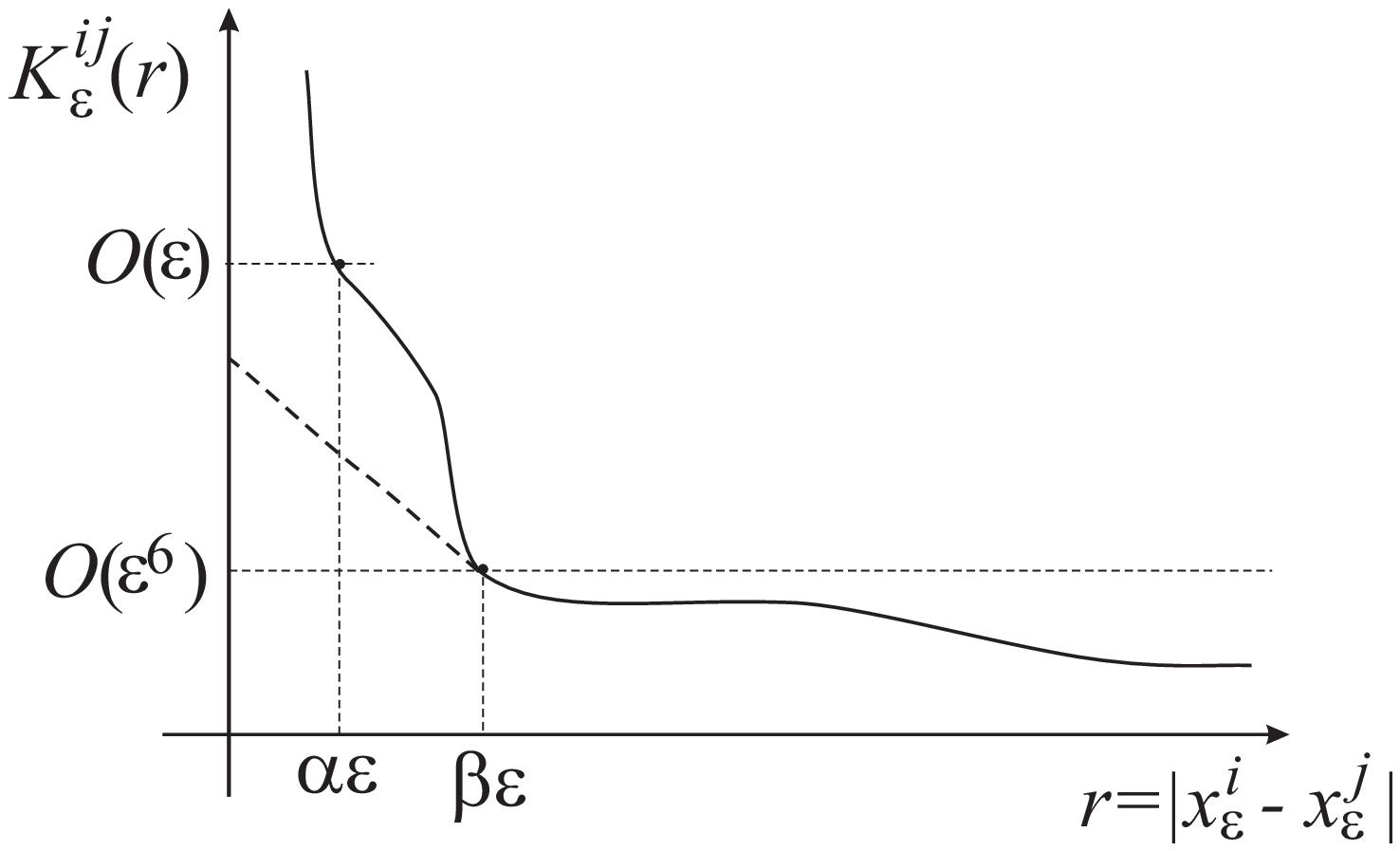}}
   \caption{}
   \label{R4:Ris1}
\end{figure}

The interaction energy of the system $ M_\varepsilon $ \eqref{1.1}
- \eqref{1.3} is invariant under rotations and shear. Therefore,
the equilibrium state $ (x_\varepsilon^1, ...,
x_\varepsilon^{N_\varepsilon})$ of the system is not isolated:
rotations and shifts are allowed. To exclude this we fix the part
of the particles $ M_\varepsilon^0 \subset M_\varepsilon $ on the
boundary $\partial \Omega$  (at the corresponding points $
x_\varepsilon^i \in \partial \Omega $ $ u_\varepsilon^i = 0 $). We
assume the following conditions hold.

I. {\underline{The condition of "$\varepsilon$-net" on the
boundary $\partial\Omega$}}.  The set $ M_\varepsilon^0 $ of
particles  assigned to $ \partial \Omega $ is a $\varepsilon$-net
for $ \partial \Omega$. It is clear that the number of such
particles is $ N_\varepsilon^0 = O (\varepsilon^{- 2}) \ll
N_\varepsilon $

II. {\underline{The triangulation condition}}.

Let $ \Gamma_\varepsilon $ be a graph with vertices at points $
x_\varepsilon^i $ and edges $ (x_\varepsilon^i, x_\varepsilon^j) $
($ i, j = 1, ..., N_\varepsilon $, $ i \not= j $). Assume that for
any $ \varepsilon> 0 $ there exists a subgraph
$\Gamma_\varepsilon^\prime \subset \Gamma_\varepsilon $ with the
same set of vertices $ M_\varepsilon $ and edges of length $ \vert
x_\varepsilon^i - x_\varepsilon^j \vert = d^{ij} \varepsilon $ ($
0 <d_1 \le d^{ij} <d_2 $), that correspond to the interaction
coefficients $ K_\varepsilon^{ij} = k^{ij} \varepsilon $ ($ 0 <a
\le k^{ij} \le A $). The subgraph $ \Gamma_\varepsilon^\prime $
triangulates the domain $\Omega $. The volumes $\vert
P_\varepsilon^\alpha \vert $ of the corresponding simplexes of the
triangulation $ P_\varepsilon^\alpha $ ($ \alpha = 1 ... \hat
N_\varepsilon $) satisfy the inequality $ \vert
P_\varepsilon^\alpha \vert> C \varepsilon^3 $ ($ C> 0 $).

Under these conditions, the equilibrium state $(x_\varepsilon ^ 1,
..., x_\varepsilon^{N_\varepsilon}) $ is isolated. In the small
neighborhood of $(x_\varepsilon ^ 1, ...,
x_\varepsilon^{N_\varepsilon}) $ the nonstationary oscillations of
the system $ M_\varepsilon $ are described by the following
problem
\begin{equation}
m_\varepsilon^i\ddot u_\varepsilon^i=-\nabla_{u_\varepsilon^i}
H_\varepsilon(u_\varepsilon^1,...,u_\varepsilon^{N_\varepsilon}),
\quad x_\varepsilon^i\in\Omega, t>0, \label{1.4}
\end{equation}
\begin{equation}
u_\varepsilon^i(t)=0,\quad x_\varepsilon^i \in \partial\Omega, t>0,
\label{1.5}
\end{equation}
\begin{equation}
u_\varepsilon^i(0)=a_\varepsilon^i,\quad \dot u_\varepsilon^i(0)=b_\varepsilon^i,\quad i=1,...,N_\varepsilon,
\label{1.6}
\end{equation}
where $m_\varepsilon^i$ is a mass of $i$-th particle,
$a_\varepsilon^i$ are the given initial displacements of the
particles, $b_\varepsilon^i$ are the given initial velocities
($a_\varepsilon^i=0$, $b_\varepsilon^i=0$ when
$x_\varepsilon^i\in\partial\Omega$). There exists a unique
solution $ \{u_\varepsilon \} = \{ u_\varepsilon^1, ...,
u_\varepsilon^{N_\varepsilon} \} $ of this problem. The main goal
of the paper is to study the asymptotic behavior of the solution
as $ \varepsilon \to 0 $. We obtain a homogenized system of
equations. This system describes the leading term of the
asymptotic and is a macroscopic model of the oscillation of an
elastic medium with a nonlocal interaction.

\section{Quantitative characteristics of the system of interacting particles and formulation of main result}
\setcounter{equation}{0}

We denote by $K_h^x = K (x, h) $ cubes with centers at points $ x
\in \Omega $ and sides of length $h$ with a fixed orientation. It
is assumed that $ 0 <\varepsilon \ll h \ll 1 $  and the cube
$K_h^x $ contains a large number  of particles (of order $ O \left
(\frac {h^3} {\varepsilon^3} \right ) $). Consider the following
functional of the symmetric tensor $ T = \{T_{np} \}_{n, p = 1}^3
$:
\begin{equation}
\begin{aligned}
H_{K_h^x}(T)=\hskip 5cm\\
=\inf\limits_{v_\varepsilon} \bigg \{ \frac{1}{2}
\sideset{}{_{K_h^x}}\sum_{{i,j}\atop \vert x_\varepsilon^i -
x_\varepsilon^j\vert \le
\beta\varepsilon}\,<E_\varepsilon^{ij}(v_\varepsilon^i-v_\varepsilon^j),
(v_\varepsilon^i-v_\varepsilon^j)> +
\\+\sideset{}{_{K_h^x}}\sum_i
h^{-2-\gamma}\left\vert
v_\varepsilon^i-\sum_{n,p=1}^3\psi^{np}(x_\varepsilon^i)T_{np}\right\vert^2\bigg
\}.
\end{aligned}
\label{2.1}
\end{equation}
The sum $ \sideset{}{_ {K_h^x}} \sum $ consist of particles
$x_\varepsilon^i \in K_h^x $ and $ \inf $ is taken over
displacements $ v_\varepsilon = \{v_\varepsilon^i, i = 1, ...,
N_\varepsilon \} $ of these particles. The vector function $
\psi^{np} (x) $ is defined by equality $\psi^{np} ( x) = \frac {1}
{2} (x_n e^p + x_p e^n) $, and $ \gamma $ is an arbitrary penalty
parameter: $ 0 <\gamma <2 $.

The functional $H_{K_h^x}(T)$ is quadratic and we can rewrite it
in the form
\begin{equation}
H_{K_h^x}(T)=\sum_{n,p,q,r=1}^3a_{npqr}(x;\varepsilon,h;\gamma)T_{np}T_{qr},
\label{2.2}
\end{equation}
where $a_{npqr}(x;\varepsilon,h;\gamma)$ are the components of the
symmetric tensor of 4-th rank in $\mathbb{R}^3$:
$a_{npqr}=a_{qrnp}=a_{pnqr}=...$. This tensor is a mesoscopic
($0<\varepsilon\ll h\ll 1$) characteristic of the concentration of
the short-range interaction energy in a neighborhood of the point
$x\in\Omega$

Assume that the limits
\begin{equation}
\lim\limits_{h\to 0}\lim\limits_{\overline{\varepsilon\to 0}}\frac{a_{npqr}(x;\varepsilon,h;\gamma)}{h^3}=
   \lim\limits_{h\to 0}\overline{\lim\limits_{\varepsilon\to 0}}\frac{a_{npqr}(x;\varepsilon,h;\gamma)}{h^3}=a_{npqr}(x)
\label{2.3}
\end{equation}
exist.

{\underline{Remark}}. Formally, the limit tensor $ \{a_{npqr} (x)
\}_ {n, p, q, r = 1}^3 $ must depend on the parameter $ \gamma $
and the orientation of the cubes $ K (x, h) $. But the main result
and the example in Section 6 show that the limiting tensor $
\{a_{npqr} (x) \}_ {n, p, q, r = 1}^3 $ does not depend on the
parameter $ \gamma $ and the orientation of the cubes $ K (x, h) $

Let $\rho_\varepsilon (x)\in L_\infty (\Omega)$ be a density of
the distribution of particles masses and let $\varphi_\varepsilon
(x,y) \in L_\infty (\Omega\times\Omega)$ be a function of the
distribution of the pairs of particles in $\Omega\times \Omega$
with long-range interaction. We will denote by $V^i_\varepsilon$
($i=1,\dots, N_\varepsilon$) Voronoi cells of a set of points
$x^i_\varepsilon \in \Omega$
$$
V^i_\varepsilon = \bigcap\limits_{j=1}^{N_\varepsilon} \{
x\in\Omega: \vert x-x^i_\varepsilon\vert <\vert
x-x^j_\varepsilon\vert\},
$$
$\vert V^i_\varepsilon\vert$ denotes the volume of the cell and
$\chi^i_\varepsilon (x)$ is a characteristic function of the cell.
Assume that
\begin{equation}
\rho_\varepsilon(x)=\sum_{i=1}^{N_\varepsilon}\frac{m_\varepsilon^i}{\vert
V_\varepsilon^i\vert} \chi_\varepsilon^i(x), \label{2.4}
\end{equation}
\begin{equation}
\varphi_\varepsilon(x,y)=\varepsilon^6\sum_{i,j=1\atop \vert
x_\varepsilon^i-x_\varepsilon^j\vert  \ge
\beta\varepsilon}^{N_\varepsilon}\frac{A_\varepsilon^{ij}}{\vert
V_\varepsilon^i\vert\vert V_\varepsilon^j\vert} \chi_\varepsilon^i
(x) \chi_\varepsilon^j (y) \label{2.5}
\end{equation}
where $m_\varepsilon^i$ are the masses of the particles,
$A_\varepsilon^{ij}$ are the elements of the adjacency matrix
$A_\varepsilon =\{A_\varepsilon^{ij}\}_{i,j=1}^{N_\varepsilon}$ of
the complete graph $\Gamma_\varepsilon$ for the system
$M_\varepsilon$ (see \eqref{1.3}).

Suppose, that for any $i=1,...,N_\varepsilon$
\begin{equation}
m_\varepsilon^i=m^i\varepsilon^3 \quad\quad (0<m_1\le
m^i_\varepsilon\le m_2<\infty). \label{2.6}
\end{equation}
By the triangulation condition II $\vert
V_\varepsilon^i\vert=c_\varepsilon^i\varepsilon^3$ ($0<C_1\le
c^i\le C_2<\infty$), and the  estimates
$\Vert\rho_\varepsilon\Vert_{L_\infty(\Omega)}<C$,
$\Vert\varphi_\varepsilon\Vert_{L_\infty(\Omega\times\Omega)}<C$
are valid uniformly with respect to $\varepsilon$. Hence the set
of functions $\{\rho_\varepsilon (x), \varepsilon
>0\}$ is *-weakly compact in $L_\infty(\Omega)$ and the set $\{\varphi_\varepsilon (x,y), \varepsilon >0\}$ is
*-weakly compact in $L_\infty (\Omega\times\Omega)$ (see
\cite{21}, \cite{22}).

We assume that
\begin{equation}
\rho_\varepsilon(x)\rightharpoonup\rho(x)\quad
{\text{*-weakly}}\quad {\text{in}}\,\,\,L_\infty(\Omega),
\label{2.7}
\end{equation}
\begin{equation}
\varphi_\varepsilon (x,y)\rightharpoonup \varphi(x,y) \quad
{\text{*-weakly}}\quad
{\text{in}}\,\,\,L_\infty(\Omega\times\Omega), \label{2.8}
\end{equation}
as $\varepsilon\to 0$. Here $\rho(x)>0$ and $\varphi(x,y)\ge 0$
are the functions in $L_\infty(\Omega)$ and
$L_\infty(\Omega\times\Omega)$ respectively.

For each discrete function
$u_\varepsilon(x)=\{u_\varepsilon^1,...,u_\varepsilon^{N_\varepsilon}\}$
that defined at the points $x_\varepsilon^i$:
$u_\varepsilon(x_\varepsilon^i)=u_\varepsilon^i$ we will match the
vector function $\tilde u_\varepsilon (x)\in L_\infty(\Omega)$ by
the formula
\begin{equation}
\tilde u_\varepsilon(x)=\sum_{i=1}^{N_\varepsilon}u_\varepsilon^i\chi_\varepsilon^i(x).
\label{2.9}
\end{equation}
The vector-functions $\tilde a_\varepsilon(x)\in
L_\infty(\Omega)$, $\tilde b_\varepsilon(x) \in L_\infty(\Omega)$
correspond to the initial data
$\{a_\varepsilon^1,...,a_\varepsilon^{N_\varepsilon}\}$ and
$\{b_\varepsilon^1,...,b_\varepsilon^{N_\varepsilon}\}$ in
\eqref{1.4}-\eqref{1.6}. The vector-function $\tilde u(x,t)\in
L_\infty(\Omega\times[0,T])$ $\forall T>0$ correspond to the
solution $\{
u_\varepsilon^1(t),...,u_\varepsilon^{N_\varepsilon}(t)\}$ of the
problem.

We assume that
\begin{equation}
\tilde a_\varepsilon(x)\to a(x),\quad \tilde b_\varepsilon (x)\to
b(x) \quad  {\textrm{in}}\,\,\,L_2(\Omega), \label{2.10}
\end{equation}
as $\varepsilon\to 0$. Here $a(x)$ and $b(x)$ are the vector
functions from $\stackrel{\circ}{W}_2^1(\Omega)$. Suppose that the
inequality
\begin{equation}
\sum_{i,j=1}^{N_\varepsilon}\langle
E_\varepsilon^{ij}(a_\varepsilon^i-a_\varepsilon^j)
(a_\varepsilon^i-a_\varepsilon^j)\rangle <C \label{2.11}
\end{equation}
holds uniformly with respect to $\varepsilon$.

 Now we can formulate the main result.

\begin{theorem}
Let the system of interacting particles $M_\varepsilon $ with the
interaction energy \eqref{1.1} - \eqref{1.3} and the masses $
m_\varepsilon^i $ \eqref{2.6} be located in $\bar \Omega $ and
conditions {\rm {I}} and {\rm {II}} are fulfilled. Suppose that
conditions \eqref{2.3}, \eqref{2.7}, \eqref{2.8} and \eqref{2.10},
\eqref{2.11} hold as $\varepsilon \to 0 $. Then the vector
function $ \tilde u_\varepsilon (x, t) $ constructed by
 \eqref{2.9} using the solution $ u_\varepsilon (t) = \{u_\varepsilon^1 (t), ...,
u_\varepsilon^{N_\varepsilon} (t) \} $ of the problem \eqref {1.4}
- \eqref {1.6} converges in $L_2 (\Omega \times [0, T])$ as
$\varepsilon \to 0 $ to  the solution $ u (x, t) $ of the
following initial-boundary value problem
\begin{equation}
\begin{aligned}
\rho(x)\frac{\partial^2 u}{\partial t^2} -\sum_{n,p,q,r=1}^3 \frac{\partial}{\partial x_q}\{a_{npqr}(x) e_{np}[u]e^r\} +\\
+\int\limits_\Omega G(x,y)(u(x,t)-u(y,t))dy =0,\quad
x\in\Omega,t>0,
\end{aligned}
\label{2.12}
\end{equation}
\begin{equation}
u(x,t)=0,\quad\quad x\in\partial\Omega,t>0,
\label{2.13}
\end{equation}
\begin{equation}
u(x,0)=a(x),\quad \frac{\partial u}{\partial t}(x,0)=b(x).
\label{2.14}
\end{equation}
Here $e_{np}[u]=\frac{1}{2} \left ( \frac{\partial u_n}{\partial
x_p}+ \frac{\partial u_p}{\partial x_n}\right )$ are the
components of the elasticity tensor, $e^r$ is the unit vector of
$x_r$ axis, and the elements of the matrix $G(x,y)$ are defined by
$$
G_{kl} (x,y) = \frac{K(\vert x-y\vert) \varphi(x,y)}{\vert x
-y\vert^2} (x_k-y_k)(x_l-y_l).
$$
\end{theorem}

The proof of the theorem is carried out in Sections 4,5. In
remainder of this Section we give the main ideas of the prove. By
the Laplace transform in time we reduce in Section 4 the initial
problem to a stationary problem with a spectral parameter $
\lambda $ ($ {\rm {Re}} \lambda> 0 $). We formulate the
variational formulation of the problem for real $ \lambda > 0 $.
Then we study the asymptotic behavior of its solution as $
\varepsilon \to 0 $ and obtain the homogenized equation. Using the
Vitali's theorem we investigate in Section 5 analytic properties
of solutions of the initial and homogenized stationary problems on
$ \lambda $ for $ {\rm {Re}} \lambda> 0 $ . We prove the
convergence of the solutions and, finally, we prove the
convergence of solutions of the original non-stationary problem
\eqref{1.3} - \eqref{1.6} to the solution of the homogenized
problem \eqref{2.11} - \eqref{2.13} with the help of the inverse
Laplace transform, .

\section{Auxiliary propositions}
\setcounter{equation}{0}

Let us denote by $L_\varepsilon^i(x)$ a continuous function in
$\mathbb{R}^3$ that is linear in every simplex
$P_{k\varepsilon}^\alpha$ (condition of triangulation II) and
$L_\varepsilon^i(x_\varepsilon^j)=\delta_{ij}$ at
$x_\varepsilon^i$. It is clear that  it is non-zero only in
simplexes with vertices  $x_\varepsilon^i$.

Using this function we construct a piecewise linear spline $\hat
u_\varepsilon (x)$ to interpolate the given discrete
vector-function
$u_\varepsilon=\{u_\varepsilon^1,...,u_\varepsilon^{N_\varepsilon}\}$:
\begin{equation}
\hat u_\varepsilon(x)=\sum_{i=1}^{N_\varepsilon}u_\varepsilon^i
L_\varepsilon^i(x), \label{3.1}
\end{equation}
where $u_\varepsilon^i=u_\varepsilon(x_\varepsilon^i)$.

In what follows, we assume that $u_\varepsilon^i=0$ for
$x_\varepsilon^i\in\partial\Omega$. Then, $ \hat u_\varepsilon (x)
\in \stackrel {\circ} {W}_2^1 (\Omega) $ for any $ \varepsilon >
0$ if the domain $ \Omega $ is convex. If $\Omega $ is not convex,
then $ \hat u_\varepsilon (x) \in \stackrel {\circ} {W}_2^1
(\Omega_\delta)$ for sufficiently small $\varepsilon \le
\varepsilon (\delta)$. Here $ \Omega_\delta $ is a domain in
$\mathbb {R}^3 $ such that $ \Omega \Subset \Omega_\delta $ and
${\textrm {dist}} (\partial \Omega, \partial \Omega_\delta) =
\delta $ ($\forall \delta > 0 $). This statement follows from
conditions I, II, and smoothness of the $ \partial \Omega $.

\begin{lemma}
Let us construct vector-functions $ \tilde u_\varepsilon (x) $ and
$\hat u_\varepsilon (x) $ by formulas \eqref{2.9} and \eqref{3.1}
for the same set of vectors $( u_\varepsilon^1, ...,
u_\varepsilon^{N_ \varepsilon})$ ($ u_\varepsilon^i = 0$ for
$x_\varepsilon^i \in \partial \Omega $). If the inequality
$$
\Vert \hat u_\varepsilon\Vert_{W_2^1(\Omega)}<C,
$$
holds uniformly with respect to $ \varepsilon $, then
$$
\Vert \hat u_\varepsilon - \tilde
u_\varepsilon\Vert_{L_2(\Omega)}\to 0\quad {\textrm{ as}}\quad
\varepsilon\to 0.
$$
\end{lemma}

{\underline{Proof}}. Denote by $ v_\varepsilon (x) = \hat
u_\varepsilon (x) - \tilde u_\varepsilon (x) $. Let $
V_\varepsilon^i $ be the Voronoi cell at the point $
x_\varepsilon^i $, and $ P_\varepsilon^\alpha $ be a simplex with
vertex at the point $ x_\varepsilon^i $ (see condition of
triangulation II). By \eqref{3.1} with $x\in
V_\varepsilon^i\bigcap P_\varepsilon^\alpha$, we get
\begin{equation}
\vert\nabla v_\varepsilon(x)\vert^2=\vert\nabla \hat
u_\varepsilon(x)\vert^2\equiv\vert
M_\varepsilon^{i\alpha}\vert^2=const. \label{3.2}
\end{equation}

Taking into account $v_\varepsilon (x_\varepsilon^i)=0$, we obtain
$$
v_\varepsilon(x)=\int\limits_0^{\vert
x-x_\varepsilon^i\vert}\frac{\partial v_\varepsilon}{\partial r}
(x_\varepsilon^i+r(x-x_\varepsilon^i))dr,\quad x\in
V_\varepsilon^i\bigcap P_\varepsilon^\alpha.
$$
By this equality, condition II and \eqref{3.2}, we have
$$
\vert v_\varepsilon(x)\vert^2\le C\varepsilon^2\vert
M_\varepsilon^{i\alpha}\vert^2,  \quad x\in V_\varepsilon^i\bigcap
P_\varepsilon^\alpha
$$
and, consequently
$$
\int\limits_\Omega \vert v_\varepsilon (x)\vert^2 dx =\sum_{i,j}
\int\limits_{V_\varepsilon^i\bigcap  P_\varepsilon^\alpha} \vert
v_\varepsilon (x)\vert^2 dx\le C\varepsilon^2
\sum_{i,\alpha}\int\limits_{V_\varepsilon^i\bigcap
P_\varepsilon^\alpha} \vert M_\varepsilon^{i\alpha}\vert^2 dx.
$$
Thus, according to \eqref{3.2} the inequality
$$
\int\limits_\Omega \vert v_\varepsilon (x)\vert^2 dx \le C\varepsilon^2 \int\limits_\Omega \vert \nabla \hat u_\varepsilon\vert^2 dx,
$$
holds, which establishes the assertion of the lemma.

Consider the function $G_{\varepsilon kl} (x,y)\in L_\infty
(\Omega\times\Omega)$ ($k,l=1,2,3$)
\begin{equation}
G_{\varepsilon kl} (x,y)=\varepsilon^6\sum_{i,j=1\atop \vert
x_\varepsilon^i-x_\varepsilon^j\vert >
\beta\varepsilon}^{N_\varepsilon}\frac{K(\vert
x_\varepsilon^i-x_\varepsilon^j\vert)e_{\varepsilon
k}^{ij}e_{\varepsilon l}^{ij}}{\vert V^i_\varepsilon\vert\vert
V^j_\varepsilon\vert} A_\varepsilon^{ij} \chi_\varepsilon^i(x)
\chi_\varepsilon^j(y), \label{3.3}
\end{equation}
where $e_{\varepsilon k}^{ij}$ are $k$-th components of the
vectors $e_\varepsilon^{ij}=
(x_\varepsilon^i-x_\varepsilon^j)\vert
x_\varepsilon^i-x_\varepsilon^j\vert^{-1}$.

\begin{lemma}
Let conditions \eqref{2.8} hold, then the function $G_{\varepsilon
kl} (x,y)$ converges to the function
\begin{equation}
G_{kl} (x,y) = \frac{K(\vert x-y\vert)\varphi(x,y)}{\vert
x-y\vert^2} (x_k-y_k) (x_l-y_l). \label{3.4}
\end{equation}
*-weakly in $L_\infty (\Omega\times\Omega)$ as $\varepsilon\to 0$.
\end{lemma}

{\underline{Proof}}. Let $f(x,y)$ be an arbitrary function in $L_1
(\Omega\times\Omega)$. By \eqref{3.3}, we write
\begin{equation}
\begin{aligned}
\int\limits_\Omega \int\limits_\Omega G_{\varepsilon kl} (x,y)f(x,y)dxdy =\hskip 3cm\\
=\int\limits_\Omega \int\limits_\Omega \left (
\varepsilon^6\sum_{i,j=1\atop \vert
x_\varepsilon^i-x_\varepsilon^j\vert\ge\beta\varepsilon}
\frac{A_\varepsilon^{ij}} {\vert V_\varepsilon^i\vert \vert
V_\varepsilon^j\vert}\chi_\varepsilon^i(x)
\chi_\varepsilon^j(y)\right )
R_{kl}(x,y) f(x,y) dxdy +\\
+\int\limits_\Omega \int\limits_\Omega
\varepsilon^6\sum_{i,j=1\atop \vert
x_\varepsilon^i-x_\varepsilon^j\vert>\delta} \left ( R_{kl}(
x_\varepsilon^i,x_\varepsilon^j-R_{kl}(x,y)\right )
\frac{A_\varepsilon^{ij}}{\vert V_\varepsilon^i\vert
\vert V_\varepsilon^j\vert}\chi_\varepsilon^i(x) \chi_\varepsilon^j(y) f(x,y) dx dy +\\
+\int\limits_\Omega \int\limits_\Omega
\varepsilon^6\sum_{i,j=1\atop \beta\varepsilon < \vert
x_\varepsilon^i-x_\varepsilon^j\vert\le\delta} \left ( R_{kl}(
x_\varepsilon^i,x_\varepsilon^j-R_{kl}(x,y)\right )
\frac{A_\varepsilon^{ij}}{\vert V_\varepsilon^i\vert \vert
V_\varepsilon^j\vert}\chi_\varepsilon^i(x) \chi_\varepsilon^j(y)
f(x,y) dx dy =\\
= I_{kl}^{\varepsilon 1}+I_{kl}^{\varepsilon 2} (\delta) +
I_{kl}^{\varepsilon 3} (\delta).\hskip 5cm
\end{aligned}
\label{3.5}
\end{equation}
Here
\begin{equation}
R_{kl} (x,y) = K(x,y) \frac{x_k-y_k)(x_l-y_l)}{\vert x-y\vert^2},
\label{3.6}
\end{equation}
and $\delta$ is an arbitrary number $\delta\gg\varepsilon$.

Since $R_{kl}(x,y)f(x,y)\in L_\infty (\Omega\times\Omega)$,
\begin{equation}
\lim\limits_{\varepsilon\to 0} I_{kl}^{\varepsilon 1}
=\int\limits_\Omega \int\limits_\Omega G_{kl}(x,y) f(x,y) dx dy.
\label{3.7}
\end{equation}
This follows from \eqref{2.5}, \eqref{2.8}, and \eqref{3.4},
\eqref{3.6}

As $f(x,y)\in L_1 (\Omega\times\Omega)$, the function $R_{kl}
(x,y)$ is continuous for $\vert x-y\vert >\delta
>0$, and ${\rm{diam}} V_\varepsilon^i\le C\varepsilon$, $\vert
V^i_\varepsilon\vert \ge C_2 \varepsilon^3$ (condition II). We
have
\begin{equation}
\lim\limits_{\varepsilon\to 0} I_{kl}^{\varepsilon 2}(\delta) =0,
\label{3.8}
\end{equation}
and
\begin{equation}
\lim\limits_{\delta\to 0}\limsup\limits_{\varepsilon\to 0}
I_{kl}^{\varepsilon 2}(\delta) =0 \label{3.9}
\end{equation}
for any fixed $\delta >0$.

From \eqref{3.5}-\eqref{3.9} follows the assertion of the lemma.

The following lemma plays a fundamental role in studying the
compactness of discrete vector-valued functions. The same role
plays the well-known Korn inequality for the functions in $
\stackrel {\circ} {W}^1_2 (\Omega)$ \ cite {24}.

\begin{lemma}[discrete analogue of the Korn inequality] Let conditions {\rm{I}} and {\rm{II}} hold.
Then
$$
{\sum_{i,j}}^\prime \langle
E_\varepsilon^{ij}[u_\varepsilon^i-u_\varepsilon^j],
[u_\varepsilon^i-u_\varepsilon^j]\rangle\ge C\Vert \hat
u_\varepsilon\Vert^2_{W^1_2(\Omega)}  \ge C_1\left (
\varepsilon{\sum_{i,j}}^\prime \vert
u_\varepsilon^i-u_\varepsilon^j\vert^2+\varepsilon^3{\sum_i}^\prime
\vert u_\varepsilon^i\vert^2\right ),
$$
for any discrete function $u_\varepsilon (x)$ defined at points
$x_\varepsilon^i$ by $u_\varepsilon (x_\varepsilon^i)=
u_\varepsilon^i$ $i=1,...,N_\varepsilon$, and $u_\varepsilon^i=0$
for $x_\varepsilon^i\in\partial\Omega$. Here $E_\varepsilon^{ij}$
are the pair interaction matrices (see \eqref{1.1}-\eqref{1.3});
the sum ${\sum_{i,j}}^\prime$ is taken over all $(i,j)$ of the
edges $(x_\varepsilon^j,x_\varepsilon^j)$ of triangulation
subgraph $\Gamma_\varepsilon^\prime$, and $C$, $C_1>0$ are the
constants that do not depend on $\varepsilon$.
\end{lemma}

The proof of lemma 3.3 is carried out in \cite{19}.

Next lemma establish the estimates of the solution
$\{v_\varepsilon^i\}$ of the problem \eqref{2.1}. We give this
lemma without proof. For more details we refer to \cite{19}.

\begin{lemma}
Let conditions \eqref{2.2} hold. Then
$$
\underset{i,j\quad\atop \vert
x_\varepsilon^i-x_\varepsilon^j\vert\le\beta\varepsilon}{{\sum}_{K_h^x}}
\langle E_\varepsilon^{ij}(v_\varepsilon^i-v_\varepsilon^j),
(v_\varepsilon^i-v_\varepsilon^j)\rangle =O(h^3),
$$
$$
\underset{i\quad}{{\sum}_{K_h^x}} \left\vert
v_\varepsilon^i-\sum_{n,p}\psi^{np}(x_\varepsilon^i)T_{np}
\right\vert^2\le O(h^{5+\gamma}),
$$
$$
\underset{\vert
x_\varepsilon^i-x_\varepsilon^j\vert\le\beta\varepsilon}{{\sum}_{K_h^x\setminus
K_{h^1}^x}} \langle
E_\varepsilon^{ij}(v_\varepsilon^i-v_\varepsilon^j),
(v_\varepsilon^i-v_\varepsilon^j)\rangle =o(h),
$$
$$
\underset{i\quad\quad\quad}{{\sum}_{K_h^x\setminus
K_{h^1}^x}}\left\vert v_\varepsilon^i-\sum_{n,p}
\psi^{np}(x_\varepsilon^i)T_{np}\right\vert^2 =o(h^{5+\gamma}),
$$
where $v_\varepsilon^i$ is a solution of the problem \eqref{2.1};
$h^1=h-2h^{1+\gamma/2}$, $\gamma>0$ and $\varepsilon\le
\hat\varepsilon(h)$.
\end{lemma}

\section{Variational formulation of the problem and asymptotic
behavior of the solution\\ as $\varepsilon\to 0$}

\setcounter{equation}{0}

By Laplace transform we convert the function $u_\varepsilon^i(t)$
of a real variable $t$ to the function of a complex variable
$\lambda$:
$$
u_\varepsilon^i(t)\to
u_\varepsilon^i(\lambda)=\int\limits_0^\infty
u_\varepsilon^i(t)e^{-\lambda t}dt,\quad  i=1,...,N_\varepsilon,
{\text{Re}}\lambda>0.
$$
Applying the Laplace transform to the problem
\eqref{1.4}-\eqref{1.6} and taking into account \eqref{1.1}, we
get the stationary problem for
$u_\varepsilon(\lambda)=\{u_\varepsilon^1(\lambda), ...,
u_\varepsilon^{N_\varepsilon}(\lambda)\}$ with a spectral
parameter $\lambda$
\begin{equation}
\begin{aligned}
\lambda^2 m_\varepsilon^i u_\varepsilon^i(\lambda)+\sum_{j=1}^{N_\varepsilon}
E_\varepsilon^{ij}(u_\varepsilon^i(\lambda)-u_\varepsilon^j(\lambda))=m_\varepsilon^if_\varepsilon^i(\lambda),\quad x_\varepsilon^i\in\Omega,\\
u_\varepsilon^i(\lambda)=0,\quad x_\varepsilon^i\in\partial\Omega,
\end{aligned}
\label{4.1}
\end{equation}
where $f_\varepsilon^i(\lambda)=\lambda
a_\varepsilon^i+b_\varepsilon^i$, $i=1,...,N_\varepsilon$.

This problem has a unique solution for all $\lambda\in\mathbb{C}$,
except the finite number of the spectrum points $\lambda =\pm i
\mu_{\varepsilon k}$ ($\mu_{\varepsilon k} >0$,
$k=1,...,N_\varepsilon^\prime <N_\varepsilon$). For $\lambda = 0$,
the problem describes the equilibrium elastic system under the
action of the forces $ m_\varepsilon^i f_\varepsilon^i$.

Solution $u_\varepsilon =\{
u_\varepsilon^1,...,u_\varepsilon^{N_\varepsilon}\}$ of the
problem \eqref{4.1} for $\lambda^2\ge 0$ minimizes the functional
\begin{equation}
\Phi_\varepsilon[v_\varepsilon]=\frac{1}{2}\sum_{i,j=1}^{N_\varepsilon} \langle E_\varepsilon^{ij}
[v_\varepsilon^i-v_\varepsilon^j], [v_\varepsilon^i-v_\varepsilon^j]\rangle +\lambda^2\sum_{i=1}^{N_\varepsilon}
m_\varepsilon^i\vert v_\varepsilon^i\vert^2 -2\sum_{i=1}^{N_\varepsilon}m_\varepsilon^i\langle f_\varepsilon^i, v_\varepsilon^i\rangle
\label{4.2}
\end{equation}
in the space $\stackrel{\circ}{J}_\varepsilon$ of discrete
vector-functions
$v_\varepsilon(x)=\{v_\varepsilon^1,...,v_\varepsilon^{N_\varepsilon}\}$
that equal $0$ on $\partial\Omega$: $v_\varepsilon^i=0$ when
$x_\varepsilon^i\in\partial\Omega$. Thus, the vector-function
$u_\varepsilon=\{ u_\varepsilon^1,...,
u_\varepsilon^{N_\varepsilon}\}$ is the solution of the
minimization problem
\begin{equation}
\Phi_\varepsilon[u_\varepsilon]=\min\limits_{v_\varepsilon\in\stackrel{\circ}{J}_\varepsilon}\Phi_\varepsilon[v_\varepsilon].
\label{4.3}
\end{equation}

To describe the asymptotic behavior of $ u_\varepsilon $ as
$\varepsilon \to 0 $, we introduce in $\stackrel {\circ} {W_2^1}
(\Omega) $ the functional
\begin{equation}
\begin{aligned}
\Phi[v]=\int\limits_\Omega \sum_{n,p,q,r=1}^3 a_{npqr}(x) e_{np}[v]e_{qr}[v]dx +\lambda^2 \int\limits_\Omega \rho(x)\vert v\vert^2 dx +\\
+\frac{1}{2}\int\limits_\Omega\int\limits_\Omega \langle
G(x,y)(v(x)-v(y),(v(x)-v(y)\rangle dxdy -2\int\limits_\Omega
\rho(x)\langle f,v\rangle dx.
\end{aligned}
\label{4.4}
\end{equation}
Here $e_{np}[v]=\frac{1}{2}\left [ \frac{\partial v_n}{\partial
x_p} + \frac{\partial v_p}{\partial x_n}\right ]$,  tensor $\{
a_{npqr} (x)\}_{n,p,q,r=1}^3$ is given by \eqref{2.3}, functions
$\rho(x)$ and $\varphi(x,y)$ are defined by
 \eqref{2.7} and \eqref{2.8}, vector-function
$f(x)=\lambda a(x)+b(x)$ is given by \eqref{2.10}, and the
elements of the matrix $G(x,y)$ are defined by \eqref{3.4}.

Consider the minimization problem
\begin{equation}
\Phi[u]=\min\limits_{w\in\stackrel{\circ}{W_2^1}(\Omega)} \Phi[w].
\label{4.5}
\end{equation}

\begin{theorem}
Let conditions {\rm{I}}, {\rm{II}}, \eqref{1.2}, \eqref{1.3} hold
and let the limits \eqref{2.3},  \eqref{2.7}, \eqref{2.8} exist as
$\varepsilon\to 0$. Then the vector-function $\tilde u_\varepsilon
(x)$  constructed by \eqref{2.9} on the solution $u_\varepsilon
=\{ u_\varepsilon^1,...,u_\varepsilon^{N_\varepsilon}\}$ of the
minimization problem \eqref{4.3} converges in $L_2(\Omega)$ to the
solution of the minimization problem \eqref{4.5}
\end{theorem}

 {\underline{Proof.}} Taking into account that $\{0\}\in\stackrel{\circ}{J}_\varepsilon$
 and $\Phi_\varepsilon [0]=0$, we get the inequality
\begin{equation}
\begin{aligned}
\sum_{i,j=1}^{N_\varepsilon}\langle E_\varepsilon^{ij}[u_\varepsilon^i-u_\varepsilon^j],
[u_\varepsilon^i-u_\varepsilon^j]\rangle +2\lambda^2 \sum_{i=1}^{N_\varepsilon} m_\varepsilon^i\vert u_\varepsilon^i\vert^2 \le\\
\le 4\left \{ \sum_{i=1}^{N_\varepsilon} m_\varepsilon^i \vert f_\varepsilon^i\vert^2\right \}^{1/2}
\left \{ \sum_{i=1}^{N_\varepsilon} m_\varepsilon^i \vert u_\varepsilon^i\vert^2\right \}^{1/2}.
\end{aligned}
\label{4.6}
\end{equation}
From \eqref{2.4}, \eqref{2.7}, \eqref{2.10} and condition II it
follows that
\begin{equation}
\sum_{i=1}^{N_\varepsilon} m_\varepsilon^i \vert f_\varepsilon^i\vert^2 \le C (\vert\lambda\vert^2 +1),
\label{4.7}
\end{equation}
where $C$ does not depend on $\varepsilon$.

We construct the vector-function $\hat u_\varepsilon (x)$ by
\eqref{3.1} where $u_\varepsilon=\{u_\varepsilon^1,...,
u_\varepsilon^{N_\varepsilon}\}$ is the solution of \eqref{4.3} .
By inequalities \eqref{4.6}, \eqref{4.7} and lemma 3.3 we get
\begin{equation}
\Vert \hat u_\varepsilon\Vert_{W_2^1(\Omega)}\le C.
\label{4.8}
\end{equation}
The inequality is satisfied uniformly with respect to
$\varepsilon$.

Thus, the set of vector-functions  $\{ \hat u_\varepsilon
(x),\varepsilon >0\}$ is a weakly compact set in
$\stackrel{\circ}{W_2^1}(\Omega)$. We can extract a subsequence
$\{ \hat u_{\varepsilon_k} (x),\varepsilon_k \to 0\}$ converges to
the vector-function $u(x)\in\stackrel{\circ}{W_2^1}(\Omega)$
weakly in $\stackrel{\circ}{W_2^1}(\Omega)$ and strongly in
$L_q(\Omega)$ ($q\le 6$).

By \eqref{2.9} we construct the subsequence $\{ \tilde
u_{\varepsilon_k} (x),\varepsilon_k \to 0\}$ for the set of
vector-functions $\{u_{\varepsilon_k}^1,...,
u_{\varepsilon_k}^{N_\varepsilon}\}$. According to lemma 3.1 and
\eqref{4.8} the subsequence converges to $u(x)$ in $L_q(\Omega)$.
Let us prove that $u(x)$ minimizes \eqref{4.5}. For this purpose
we rewrite the functional $\Phi_\varepsilon$ \eqref{4.2} in the
form:
\begin{equation}
\Phi_\varepsilon[v_\varepsilon]=\Phi_{1\varepsilon}[v_\varepsilon]+\Phi_{2\varepsilon}[v_\varepsilon],
\label{4.9}
\end{equation}
where
\begin{equation}
\Phi_{1\varepsilon}[v_\varepsilon] =\frac{1}{2}\sum_{i,j=1\atop
\vert x_\varepsilon^i-x_\varepsilon^j\vert \le
\beta\varepsilon}^{N_\varepsilon}\langle
E_\varepsilon^{ij}(v_\varepsilon^i-v_\varepsilon^j),
(v_\varepsilon^i-v_\varepsilon^j)\rangle, \label{4.10}
\end{equation}

\begin{equation}
\begin{aligned}
\Phi_{2\varepsilon}[v_\varepsilon]=\hskip6cm\\
=\frac{1}{2}\sum_{i,j=1\atop \vert
x_\varepsilon^i-x_\varepsilon^j\vert \ge
\beta\varepsilon}^{N_\varepsilon}\langle
E_\varepsilon^{ij}(v_\varepsilon^i-v_\varepsilon^j),
(v_\varepsilon^i-v_\varepsilon^j)\rangle
+\lambda^2\sum_{i=1}^{N_\varepsilon} m_\varepsilon^i \vert
v_\varepsilon^i\vert^2-\\
-2\sum_{i=1}^{N_\varepsilon}m_\varepsilon^i\langle
f_\varepsilon^i,v_\varepsilon^i\rangle.\hskip4cm
\end{aligned}
\label{4.11}
\end{equation}
Strong interactions between nearby particles are included in
functional $\Phi_{1\varepsilon}$. According to \eqref{1.2},
\eqref{1.3} this interaction has the order $O(\varepsilon)$ (see
\cite{19}). Recalling that $\hat u_\varepsilon\to u$ converges
weakly in $\stackrel{\circ}{W_2^1}(\Omega)$ as
$\varepsilon=\varepsilon_k\to 0$, and taking into account
\eqref{2.2}, we get the lower bound for $\Phi_{1\varepsilon}$ by
the same method as in \cite{19}:
\begin{equation}
\lim\limits_{\overline{\varepsilon=\varepsilon_k\to
0}}\Phi_{1\varepsilon}[\tilde u_\varepsilon]\ge\Phi_1[u] =
\int\limits_\Omega \sum_{n,p,q,r=1}^3 a_{npqr} (x)
e_{np}[u]e_{qr}[u]dx. \label{4.12}
\end{equation}

By \eqref{1.2}, \eqref{1.3},  \eqref{2.4}, \eqref{3.5} we can
rewrite $\Phi_{2\varepsilon}[v_\varepsilon]$ in the form:
\begin{equation}
\begin{aligned}
\Phi_{2\varepsilon}[\tilde u_\varepsilon]=\int\limits_\Omega
\int\limits_\Omega \sum_{k,l=1}^3 G_{\varepsilon kl} (\tilde
u_{\varepsilon k} (x)- \tilde u_{\varepsilon k} (y)) (\tilde
u_{\varepsilon l} (x)- \tilde u_{\varepsilon l} (y))
dxdy +\\
+\lambda^2 \int\limits_\Omega \rho_\varepsilon(x)\vert \tilde
u_\varepsilon\vert^2 dx-2 \int\limits_\Omega
\rho_\varepsilon(x)\langle\tilde f_\varepsilon, \tilde
u_\varepsilon\rangle dx,
\end{aligned}
\label{4.13}
\end{equation}
where
$$
\tilde f_\varepsilon (x)=\sum\limits_{i=1}^{N_\varepsilon}
f_\varepsilon^i\chi_\varepsilon^i (x)=
\sum\limits_{i=1}^{N_\varepsilon} (\lambda a_\varepsilon^i +
b_\varepsilon^i) \chi_\varepsilon^i (x),
$$
and $\tilde u_{\varepsilon k} (x)$ is $k$-th component of the
vector-function $\tilde u_\varepsilon(x)$ \eqref{2.4}.

Since $\tilde u_\varepsilon \to u$ in $L_2(\Omega)$ as
$\varepsilon=\varepsilon_k\to 0$, we have
$$
(\tilde u_{\varepsilon k} (x)- \tilde u_{\varepsilon k} (y))
(\tilde u_{\varepsilon l} (x)- \tilde u_{\varepsilon l} (y)) \to
(u_{ k} (x)- u_{ k} (y)) (u_{ l} (x)- u_{ l} (y))\quad
{\textrm{in}}\,\, L_2 (\Omega\times\Omega)
$$
and by \eqref{2.10}
$$
\langle \tilde f_\varepsilon, \tilde u_\varepsilon \rangle \to
\langle f,u\rangle \quad {\textrm{in}}\,\, L_1 (\Omega),
$$
where $f(x)=\lambda a(x) + b(x)$.

From the above, by lemma 3.2, \eqref{2.7}, and \eqref{4.13}, we
obtain
\begin{equation}
\begin{aligned}
\lim\limits_{\varepsilon=\varepsilon_k\to 0}\Phi_{2\varepsilon}
[\tilde u_\varepsilon]= \Phi_2[u]=\hskip3cm \\
=\int\limits_\Omega \int\limits_\Omega \langle G(x,y)
(u(x)-u(y)),(u(x)-u(y))\rangle dx dy +\\
+\lambda^2 \int\limits_\Omega \rho(x) u^2 (x) dx -
2\int\limits_\Omega \rho (x) \langle f(x), u(x)\rangle dx.
\end{aligned}
\label{4.14}
\end{equation}

On account of \eqref{4.9}, \eqref{4.13}, \eqref{4.14}, we get the
lower bound for $\Phi_\varepsilon [\tilde u_\varepsilon]$:
\begin{equation}
\lim\limits_{\overline{\varepsilon=\varepsilon_k\to
0}}\Phi_\varepsilon [\tilde u_\varepsilon]\ge \Phi[u],
\label{4.15}
\end{equation}
where $u(x)$ is a limit in $L_2 (\Omega)$ of the vector-functions
$\tilde u_\varepsilon (x)$, and the functional
$\Phi[u]=\Phi_1[u]+\Phi_2[u]$ is defined by \eqref{4.4}.

In order to get the upper bound, we introduce the test
vector-function $w_{\varepsilon k}=(w_{\varepsilon
k}^1,...,w_{\varepsilon k}^{N_\varepsilon})$
 in $\stackrel{\circ}{J}_\varepsilon$ for the problem \eqref{4.3}.
 To this end, we cover $\Omega$ by cubes $K_h^\alpha
= K(x^\alpha,h)$ with centers at the points $x^\alpha$ and sides
of length $h$. The centers of the cubes form a cubic lattice with
period $h-h^{1+\gamma/2}$ ($0<\gamma <2$). By this covering we
construct the partition of unity $\varphi_\alpha (x)$. Namely, the
set of functions with the following properties: $\varphi_\alpha
(x)\in C_0^2 (K_h^\alpha)$, $\sum_\alpha \varphi_\alpha (x)=1$,
$\varphi_\alpha (x)=0$ when $x\not\in K_h^\alpha$, $\varphi_\alpha
(x)=1$ when $x\in K_h^\alpha\setminus
\bigcup\limits_{\beta\not=\alpha}K_h^\beta$; $\vert\nabla
\varphi_\alpha(x)\vert\le Ch^{-1-\gamma/2}$.

Let $w(x)$ be an arbitrary vector-functions in $C^2(\Omega)$ with
the compact support in $\Omega$. Define
\begin{equation}
\begin{aligned}
w_{\varepsilon h}^i = \sum_\alpha \bigg \{ w(x^\alpha) +\sum_{n,p=1}^3 (e_{np}[w(x^\alpha)]
v_{\varepsilon h}^{\alpha np} (x_\varepsilon^i)+\quad\quad\quad\\
 + \omega_{np}[w(x^\alpha)]\varphi^{np}(x_\varepsilon^i-x^\alpha))\bigg\}\varphi_\alpha (x_\varepsilon^i),\quad i=1,...,N_\varepsilon.
\end{aligned}
\label{4.16}
\end{equation}
Here $v_{\varepsilon h}^{\alpha np}(x_\varepsilon^i)$ is a
minimizer of the functional \eqref{2.1}  in the cube $K_h^\alpha$
for $T=T^{np}$ ($T^{np}_{ik}=\delta_{in}\delta_{pk}$);
$e_{np}[w]$, $\omega_{np}[w]$ are a symmetric and antisymmetric
parts of the tensor $\nabla w$; $\varphi^{np} (x)=\frac{1}{2}(x_n
e^p-x_pe^n)$.

Using the properties of the discrete vector-functions
$v_{\varepsilon h}^{\alpha np}$ (see lemma 3.4), the properties of
the partition of unity $\{\varphi_\alpha (x)\}$, and by
\eqref{2.2} we get
\begin{equation}
\lim\limits_{h\to 0}{\overline{\lim\limits_{\varepsilon\to 0}}} \Phi_{1\varepsilon} [w_{\varepsilon h}]\le\Phi_1[w],
\label{4.17}
\end{equation}
in the same manner as in \cite{19}. The functional $\Phi_1$ is
defined by \eqref{4.12}.

To estimate $\Phi_{2\varepsilon}[w_{\varepsilon h}]$, we use the
following equality for the vector-functions $w(x)\in
C_0^2(\Omega)$ for $x\in K_h^\alpha$
$$
w(x)=w(x^\alpha) +\sum_{n,p} e_{np}[w(x^\alpha)]\psi^{np} (x-x^\alpha)+\omega_{np}[w(x^\alpha)]\varphi^{np}(x-x^\alpha)+O(h^2).
$$

Substituting this equality in \eqref{4.16} and applying lemma 3.4,
we conclude
\begin{equation}
\lim\limits_{h\to 0}{\overline{\lim\limits_{\varepsilon\to 0}}}
\Vert w_{\varepsilon h}-w\Vert^2_{L_2(\Omega)}=0. \label{4.18}
\end{equation}

Taking into account convergence \eqref{2.7}, \eqref{2.10}, and
lemma 3.2, we get
\begin{equation}
\lim\limits_{h\to 0}{\overline{\lim\limits_{\varepsilon\to 0}}}\Phi_{2\varepsilon}[w_{\varepsilon h}]=\Phi_2[w],
\label{4.19}
\end{equation}
where $\Phi_2$ is defined by \eqref{4.14}.

Thus, by \eqref{4.9}, \eqref{4.17}, \eqref{4.19}
$$
\lim\limits_{h\to 0}{\overline{\lim\limits_{\varepsilon\to 0}}}\Phi_\varepsilon [w_{\varepsilon h}]\le \Phi[w].
$$
Recalling that $u_\varepsilon$ is the minimizer of
$\Phi_\varepsilon$ in $\stackrel{\circ}{J}_\varepsilon$ for
sufficient small $h$ ($\varepsilon <\hat\varepsilon(h)$) and
$w_{\varepsilon h}\in \stackrel{\circ}{J}_\varepsilon$ we can
write
\begin{equation}
\lim\limits_{h\to 0}{\overline{\lim\limits_{\varepsilon\to 0}}}\Phi_\varepsilon[u_\varepsilon]\le\Phi[w].
\label{4.20}
\end{equation}
Combining \eqref{4.15} and \eqref{4.20} we obtain
$$
\Phi[u]\le\Phi[w],\quad w\in C_0^2(\Omega).
$$
The inequality is valid for any vector-function
$w\in\stackrel{\circ}{W^2_1}(\Omega)$ due to the continuity of the
functional $\Phi[w]$ in $\stackrel{\circ}{W^2_1}(\Omega)$. Thus,
the limit $u(x)$ of the vector-functions $\tilde u_\varepsilon
(x)$ by subsequence $\varepsilon=\varepsilon_k\to 0$ is a solution
of minimizing problem \eqref{4.5}. Consequently, $u(x)$ is a weak
solution of the following boundary value problem
\begin{equation}
\begin{aligned}
\sum_{n,p,q,r=1}^3 \frac{\partial}{\partial x_q} \left \{a_{npqr}(x)e_{np}[u]e^r\right\} +\lambda^2\rho(x) u +\\
+\int\limits_\Omega \langle G(x,y)(u(x) - u(y))\rangle dy =\lambda a(x)+b(x), \quad x\in\Omega,\\
u(x)=0,\quad x\in\partial\Omega.
\end{aligned}
\label{4.21}
\end{equation}
Since $\lambda\ge 0$, function $\rho(x)$ and matrix-function
$G(x,y)$ are non-negative, and tensor
$\{a_{npqr}(x)\}_{n,p,q,r=1}^3$ is positive definite, the problem
have a unique solution. Thus, theorem 4.1 is proved.

\section{The convergence of solutions of the problem \eqref{4.1}
to the solution of the problem \eqref{4.21} for complex $\lambda $}
\setcounter{equation}{0}

 1. Consider the problem \eqref{4.1} for complex $\lambda$ in
 semiaxis ${\textrm{Re}}\lambda >0$. Denote by $L_\varepsilon$ a
 Hilbert space of finite sets of $N_\varepsilon$ 3-components
 complex vectors defined in $x_\varepsilon^i\in\bar\Omega$: $u_\varepsilon =\{ u_\varepsilon^1,...,
 u_\varepsilon^{N_\varepsilon}\}$. If
 $x_\varepsilon^i\in\partial\Omega$ then $u_\varepsilon^i=0$.
 Define a scalar product in $L_\varepsilon$:
 $$
 (u_\varepsilon, w_\varepsilon)_\varepsilon =\sum_{i=1}^{N_\varepsilon} \langle u_\varepsilon^i, \bar w_\varepsilon^i\rangle m_\varepsilon^i,
 $$
 where $m_\varepsilon^i$ is a mass of the point $x_\varepsilon^i$.
 By parenthesis $\langle\, ,\,\rangle$ we denote the scalar product
 in $\mathbb{R}^3$. The bar denotes the complex conjugation. The
 corresponding norm is denoted by $\Vert u_\varepsilon\Vert_\varepsilon = (u_\varepsilon, \bar u_\varepsilon )_\varepsilon^{1/2}$

Consider in $L_\varepsilon$ a linear operator $A_\varepsilon:
L_\varepsilon\to L_\varepsilon$:
 \begin{equation}
 (A_\varepsilon u_\varepsilon)_i=\begin{cases} \frac{1}{m_\varepsilon^i}\sum_{i=1}^{N_\varepsilon}E_\varepsilon^{ij}
 (u_\varepsilon^i -u_\varepsilon^j), & x_\varepsilon^i\in\Omega,\cr 0, & x_\varepsilon^i\in\partial\Omega.\end{cases}
 \label{5.1}
 \end{equation}

From \eqref{2.5}, \eqref{1.2}, \eqref{1.3} it follows that
$A_\varepsilon$ is a bounded selfadjoint operator in
$L_\varepsilon$. By lemma 3.3 $A_\varepsilon$ is positive definite
operator uniformly  with respect to $\varepsilon$:
\begin{equation}
(A_\varepsilon u_\varepsilon, u_\varepsilon)_\varepsilon
=\sum_{i,j=1}^{N_\varepsilon} \langle E_\varepsilon^{ij}
(u_\varepsilon^i-u_\varepsilon^j),
u_\varepsilon^i-u_\varepsilon^j\rangle\ge \alpha\Vert
u_\varepsilon\Vert_\varepsilon^2,\quad (\alpha >0). \label{5.2}
\end{equation}

Let us rewrite the problem \eqref{4.1} in operator form in
$L_\varepsilon$:
\begin{equation}
A_\varepsilon u_\varepsilon +\lambda^2 u_\varepsilon =\lambda a_\varepsilon +b_\varepsilon.
\label{5.3}
\end{equation}
By the indicated properties of operator $A_\varepsilon$ its
resolvent is a meromorphic operator function of the parameter
$\tau =\lambda^2$ with poles on the negative semiaxis $\tau <0$.
Hence the solution $u_\varepsilon = u_\varepsilon (\lambda)$ of
\eqref{5.3} is a holomorphic function of $\lambda$ in the
half-plane ${\textrm{Re}} \lambda >0$. Multiplying \eqref{5.3} on
$\bar u_\varepsilon$ and separating the imaginary and real parts,
taking into account \eqref{2.5} and \eqref{2.10}, we obtain the
estimate for $u_\varepsilon$ in half-plane ${\textrm{Re}} \lambda
>\sigma$ ($\forall \sigma >0$) uniform with respect to
$\varepsilon$: $\Vert u_\varepsilon\Vert_\varepsilon \le C$
($C=C(\sigma) <\infty$). This implies that the vector-function
$\tilde u_\varepsilon =u_\varepsilon (x,\lambda)$ defined by
\eqref{2.9} is a holomorphic function in ${\textrm{Re}} \lambda
>\sigma$ ($\forall \sigma >0$). Moreover, $\tilde u_\varepsilon$
is bounded in the norm of $L_2(\Omega)$ uniformly with respect to
$\varepsilon$:
\begin{equation}
\Vert \tilde u_\varepsilon \Vert_{L_2(\Omega)}\le C<\infty.
\label{5.4}
\end{equation}

2. We now turn to the problem \eqref{4.21}. Denote by $L_2(\Omega,
\rho)$ a Hilbert space of a complex-valued vector-functions in
$L_2(\Omega)$ with a weight $\rho(x) >0$. The scalar product in
$L_2(\Omega, \rho)$ we define by
$$
(u,w)_\rho =\int\limits_\Omega u(x){\overline{w(x)}} \rho(x) dx.
$$
Consider a sesquilinear form defined on the set of vector-valued
functions $C_0 (\Omega)$ that is dense in $L_2 (\Omega, \rho)$
$$
\hat A(u,w)=\frac{1}{\rho} \int\limits_\Omega \sum_{n,p,q,r=1}^3
a_{npqr} e_{np}[u]e_{qr}[\bar w]dx +
\frac{1}{2\rho}\int\limits_\Omega \langle G(x,y)[u(x)-u(y)],[\bar
w(x)-\bar w(y)]\rangle dxdy.
$$
The form generates in $L_2(\Omega, \rho)$ a selfadjoint operator
$A$, that is due to the properties of the elasticity tensor
$\{a_{npqr}\}_{n,p,q,r=1}^3$ and the long-range matrix $G(x,y)$
\cite{21}. The equality
$$
(Au,u)_\rho =\int\limits_\Omega \sum_{n,p,q,r=1}^3 a_{npqr}
(x)\vert e_{np}[u]\vert^2 dx + \frac{1}{2} \int\limits_\Omega
\langle G(x,y) [u(x)-u(y)], [\bar u(x)-\bar u(y)]\rangle dxdy,
$$
is valid. From the Korn's inequality it follows that
\begin{equation}
(Au,u)_\rho \ge C\Vert u\Vert^2_{\stackrel{\circ}{W^1_2}(\Omega)}\quad (C>0).
\label{5.5}
\end{equation}
This inequality implies that the operator $A$ is positive definite
and has a completely continuous inverse operator. Now we can
rewrite the problem \eqref{4.21} in operator form:
\begin{equation}
Au+\lambda^2u=\lambda a+b.
\label{5.6}
\end{equation}
The properties of the operator $A$ implies that equation
\eqref{5.6} has a solution $u(x)$ for complex $\lambda$
(${\textrm{Re}}\lambda >0$) and this solution is a holomorphic
function of $\lambda$ satisfying the inequality
$$
\Vert u\Vert_\rho < C.
$$

3. By theorem 2.1 the vector-function $\tilde u_\varepsilon
(x,\lambda)$ converges in $L_2(\Omega)$ for $\lambda >0$ to the
solution $u(x,\lambda)$ of the problem \eqref{4.21} (or equation
\eqref{5.6}) as $\varepsilon\to 0$. Moreover, the set of
vector-functions $\{\tilde u_\varepsilon, \varepsilon
>0\}$ is bounded by norm in $L_2(\Omega)$ uniformly with respect
to $\varepsilon$ in the half-plane ${\textrm{Re}} \lambda >\sigma$
($\forall\sigma >0$). Therefore, using Vitali's theorem and taking
into account that $u(x,\lambda)$ is holomorphic, we get the
following assertion.

\begin{theorem}
Let assumptions {\rm{I}}, {\rm{II}}, \eqref{2.2}, \eqref{2.6},
\eqref{2.9} hold. Let construct the function $\tilde u_\varepsilon
(x,\lambda)$ by \eqref{2.9} on the solution of the problem
\eqref{4.1}. Then vector-function $\tilde u_\varepsilon
(x,\lambda)$ converges in $L_2(\Omega)$ to the
solution$u(x,\lambda)$ of equation \eqref{5.6} (or the problem
\eqref{4.21}) uniformly with respect to complex $\lambda$ from the
half-plane ${\textrm{Re}} \lambda
>\sigma$ ($\forall\sigma >0$).

\end{theorem}

\section{The end of the proof of the main theorem}

\setcounter{equation}{0}

By the definition \eqref{5.1} of operator $A_\varepsilon$ the
problem \eqref{1.4} -- \eqref{1.6} is representable in
$L_\varepsilon$ in the form:
\begin{equation}
\begin{aligned}
\ddot u_\varepsilon +A_\varepsilon u_\varepsilon=0,\\
u_\varepsilon (0)=a_\varepsilon,\quad \dot u_\varepsilon (0)=b_\varepsilon.
\end{aligned}
\label{6.1}
\end{equation}
From this on account of \eqref{5.1}, we have
$$
\Vert \dot u_\varepsilon\Vert_\varepsilon^2
+\sum_{i,j=1}^{N_\varepsilon}  \langle
E_\varepsilon^{ij}(u_\varepsilon^i -u_\varepsilon^j),
(u_\varepsilon^i -u_\varepsilon^j)\rangle = \Vert
b_\varepsilon\Vert_\varepsilon^2 +\sum_{i,j=1}^{N_\varepsilon}
\langle E_\varepsilon^{ij} (a_\varepsilon^i -a_\varepsilon^j),
(a_\varepsilon^i -a_\varepsilon^j)\rangle.
$$
The equality above with the discrete Korn's inequality, the
properties of $E_\varepsilon^{ij}$ and $m_\varepsilon^i$, and
\eqref{2.9}, \eqref{2.10} implies inequality:
$$
\int\limits_{\Omega_T} \left \{ \left ( \frac{\partial \hat
u_\varepsilon}{\partial t}\right )^2  +\vert \nabla \hat
u_\varepsilon \vert^2\right \} dxdt \le CT\quad (\forall T >0),
$$
where $\hat u_\varepsilon = \hat u_\varepsilon (x,t)$ is a spline
vector-function, defined by \eqref{3.1}, and $C$ does not depend
on $\varepsilon$.

Thus the set of vector-functions $\{\hat u_\varepsilon,
\varepsilon\to 0\}$ is bounded in $W_2^1 (\Omega_T)$ uniformly
with respect to $\varepsilon$. We can extract a subsequence
$\{\hat u_\varepsilon, \varepsilon=\varepsilon_k\to 0\}$ converges
weakly in $W_2^1(\Omega_T)$ to a function $u(x,t)\in
W_2^1(\Omega_T)$ (and by embedding theorem converges strongly in
$L_q (\Omega\times (0,T))$ ($q\le 4$) and for almost all $t\in
(0,T]$ converges strongly in $L_2 (\Omega)$). By the above and
lemma 3.1 we conclude that the piecewise-constant vector-functions
$\tilde u_\varepsilon (x,t)$ defined by \eqref{2.9} converges in
$L_4 (\Omega)$ and $L_2(\Omega)$ to $u(x,t)$ for almost all $t\in
[0,T]$ as $\varepsilon=\varepsilon_k\to 0$.

Let us prove that the function $u(x,t)$ is a solution of the
problem \eqref{2.11} -- \eqref{2.13}. By definition of operator
$A$ this problem can be written in operator form:
\begin{equation}
\begin{aligned}
\ddot u+ Au=0\\
u(0)=a,\quad \dot u=b.
\end{aligned}
\label{6.2}
\end{equation}

The solution $u_\varepsilon (x,t)$ of the problem
\eqref{1.4}-\eqref{1.6} is an inverse Laplace transform of the
solution $u_\varepsilon (x,\lambda)$ of the problem \eqref{4.1}
$$
u_\varepsilon (x,t)=\frac{1}{2\pi i}\int\limits_{\sigma
-i\infty}^{\sigma +i\infty}  e^{\lambda t}u_\varepsilon
(x,\lambda) d\lambda, \quad \sigma >0.
$$
Thus we have
\begin{equation}
\tilde u_\varepsilon (x,t)=\frac{1}{2\pi i} \int\limits_{\sigma
-i\infty}^{\sigma +i\infty}  e^{\lambda t}\tilde u_\varepsilon
(x,\lambda) d\lambda, \label{6.3}
\end{equation}
where $\tilde u(x,t)$ and $\tilde u(x,\lambda)$ are defined by
\eqref{2.9}. Multiply the equality above by $\psi(x)\varphi(t)$,
where $\psi (x)\in L_2 (\Omega)$, $\varphi (t)\in C_0^2(0,T]$, and
integrate over $\Omega_T$. Changing the integration order and
integrating on $t$ by parts, we obtain
$$
\int\limits_{\Omega_T} \tilde u_\varepsilon (x,t)\psi(x)\varphi(t)
dx  dt = \frac{1}{2\pi i} \int\limits_{\sigma -i\infty}^{\sigma
+i\infty} \frac{e^{\lambda t}}{\lambda^2} \left (
\int\limits_{\Omega_T} \tilde u_\varepsilon (x,\lambda)\psi(x)
\frac{\partial^2 \varphi}{\partial t^2} dxdt\right ) d\lambda.
$$
Note that the integrals on $\lambda$ in the right-hand side
converge absolutely, due to \eqref{5.4}.

Let us pass to the limit in the equation above as
$\varepsilon=\varepsilon_k\to 0$. Passing to the limit we take
into account that $\tilde u_{\varepsilon_k} (x,t)$ converges to
$u(x,t)$ in $L_2 (\Omega)$, and $\tilde u_{\varepsilon_k}
(x,\lambda)$ converges to the solution $u(x,\lambda)$ of the
equation \eqref{5.6} in $L_2(\Omega)$ uniformly on the compacts
$\Lambda$ from the half-plane ${\text{Re}}\lambda >0$ (see theorem
2). We get
$$
\int\limits_{\Omega_T} u(x,t)\psi(x)\varphi(t) dx dt =
\int\limits_{\Omega_T}  \left \{ \frac{1}{2\pi i}
\int\limits_{\sigma -i\infty}^{\sigma +i\infty} u(x,\lambda
e^{\lambda t} d\lambda\right \} \psi(x)\varphi(t) dx dt.
$$
Since the linear combination of the functions $\psi(x)\varphi(t)$
form a dense set in $L_2(\Omega_T)$ then $u(x,t)$ is a solution of
the problem \eqref{6.2}. By the properties of operator $A$, this
problem has a unique solution. Thus $\tilde u(x,t)$ converges to
$u(x,t)$ in $L_2 (\Omega_T)$ as $\varepsilon\to 0$. The theorem
1.1 is proved.

\section{Periodic structure}
\setcounter{equation}{0}

We now consider the concrete case when the conditions of theorem
1.1 are satisfied and the elastic tensor $ \{a_{npqr} (x) \} $ and
the matrix-function $ G (x, y) $ are computed explicitly.

Suppose that the points $ x_\varepsilon^i$ of the equilibrium
state of the system are located periodically. They form a cubic
lattice with a period $ \varepsilon $. Each point
$x_\varepsilon^i$ interacts with the tops of the cube
$x_\varepsilon^j$. The points $x_\varepsilon^i$ $x_\varepsilon^j$
belongs to the same cube. For clarity, we can assume that the
interaction is carried out by elastic springs . The stiffness of
the springs (the elasticity coefficient in Hooke's law) directed
along the edges of the cubes is $ k_1 \varepsilon^2 $; directed
along the diagonals of the faces of the cubes is $ k_2 \varepsilon
$; and directed along the diagonals of the cubes is $ k_3
\varepsilon^2 $ (Fig. 2). This is a strong short-range
interaction.

\begin{figure}[h]
  \center{\includegraphics[width=0.5\linewidth]{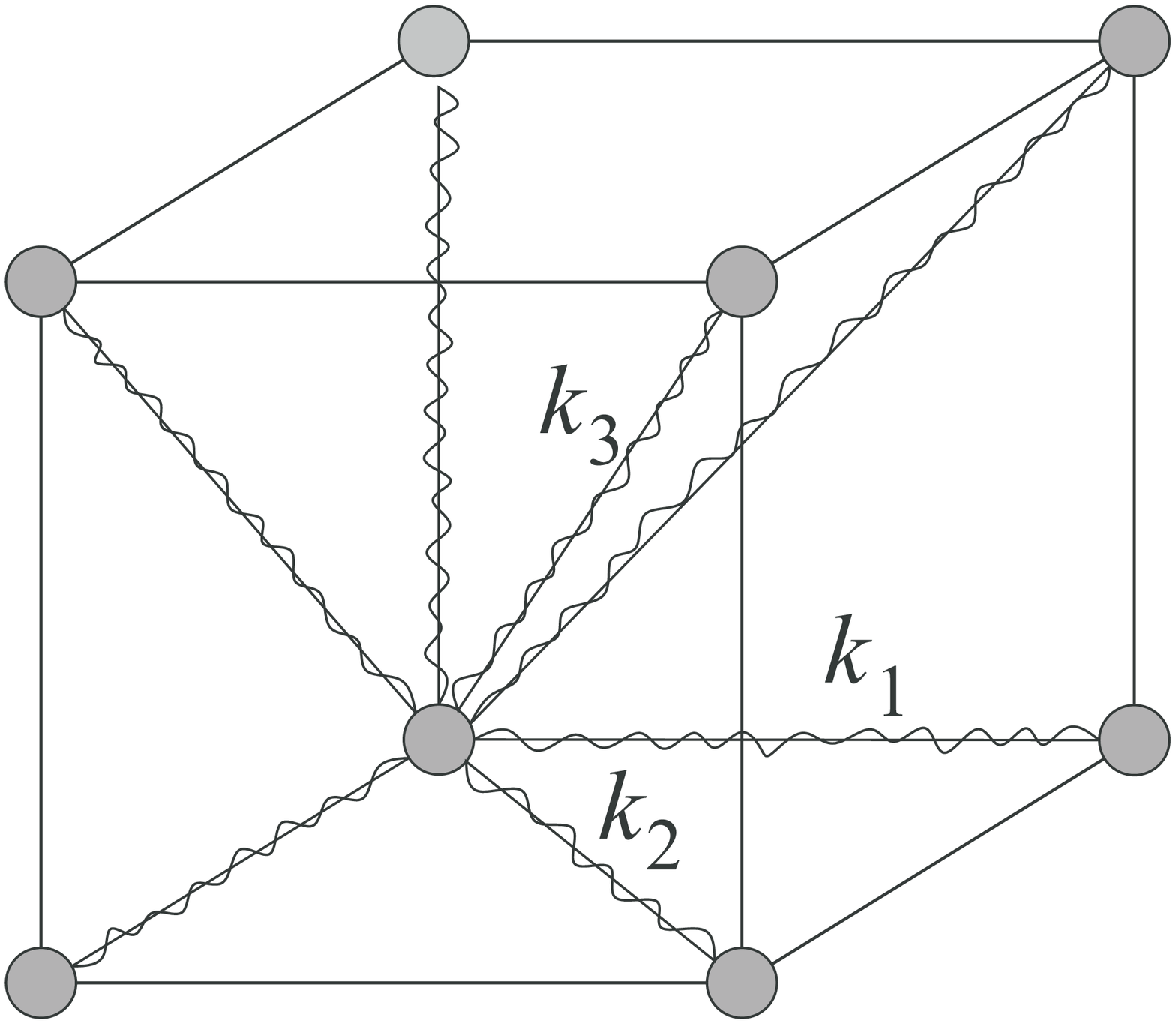}}
   \caption{}
   \label{R4:Ris2}
\end{figure}

The corresponding coefficients of interaction \eqref{1.2}
$K_\varepsilon^{ij}$ have an order $O(\varepsilon)$ $
K_\varepsilon^{ij}= k_1\varepsilon$, $4k_2\varepsilon$,
$9k_3\varepsilon$.

Let us assume that there exist a long-range interaction. Each
point $x_\varepsilon^i$ interacts with the points
$x_\varepsilon^{j}$ of the cubic sublattice
$\{x_\varepsilon^j\}^{(i)}$ with the period $N\varepsilon$
($\exists N\in \mathbb{Z}$, $N\ge 2$). This is a weak interaction
and
$$
K_\varepsilon^{ij}=\varepsilon^6 K\vert
x_\varepsilon^i-x_\varepsilon^j\vert,
$$
where $K(r)$ is a non-negative function (see \eqref{1.3}).

The system of the points $x_\varepsilon^i$ satisfies the
triangulation condition II. The corresponding interaction is
described by \eqref{1.2}, where $\alpha = \sqrt{3}$, $\beta = 2$
$K_{ij} = k_1, k_2, k_3$; $A_{ij}=1$ only for $\vert
x_\varepsilon^i-x_\varepsilon^j\vert=\varepsilon,
\sqrt{2}\varepsilon, \sqrt{3}\varepsilon$ and for $\vert
x_\varepsilon^i-x_\varepsilon^j\vert=
\sqrt{k^2+l^2+m^2}N\varepsilon$ ($k,l,m=1,2,3,\dots$).

By \eqref{2.3} the limit dense $\varphi(x,y)$ is equal to
$\frac{1}{N^3}$. Therefore, by \eqref{3.4}
\begin{equation}
G_{kl}(x,y)=\frac{K(\vert x-y\vert)(x_k-y_k)(x_l-y_l)}{N^3\vert
x-y\vert^2}. \label{7.1}
\end{equation}
The components of elasticity tensor for this system are calculated
in  \cite{19}. They are determined by formulas:
$$
a_{nnnn}= k_1+2\frac{k_2}{\sqrt{2}}+\frac{4k_3}{3\sqrt{3}}, \quad
a_{nnpp}=a_{npnp}=\frac{k_2}{\sqrt{2}}+\frac{4k_3}{3\sqrt{3}}\quad
(n\not= p)
$$
and $a_{npqr}=0$ in other cases.

{\it Remark}. If we take
$k_1=\frac{k_2}{\sqrt{2}}+\frac{8k_3}{3\sqrt{32}}$,  then the
components of limiting elasticity tensor are satisfying the
condition: $a_{nnnn}=2a_{npnp}+a_{nnpp}$ and the limit model of
elastic system is isotropic. The equation \eqref{2.12} has a form:
$$
\frac{\partial^2 u}{\partial t^2}- a\Delta u+ 2a\nabla{\rm{div}} u
+\int\limits_\Omega G(x,y) (u(x)-u(y))dy =0,
$$
where $a=a_{nnpp}=a_{npnp}$, and the elements of matrix $G(x,y)$
are defined by \eqref{7.1}.

\end{document}